\documentclass[letterpaper]{IEEEtran}

\usepackage{cite}
\usepackage{url}
\usepackage[hidelinks]{hyperref}
\usepackage{amsmath,amssymb,amsfonts}
\usepackage{graphicx}
\usepackage{array}
\usepackage{booktabs}
\usepackage{threeparttable}
\usepackage{multirow}
\usepackage{xcolor}

\providecommand{\makecell}[2][c]{\begin{tabular}[c]{@{}#1@{}}#2\end{tabular}}

\makeatletter
\def\section{\@startsection{section}{1}{\z@}%
  {1.1ex plus 0.4ex minus 0.2ex}%
  {0.5ex plus 0.3ex minus 0ex}%
  {\normalfont\normalsize\centering\scshape}}%
\def\subsection{\@startsection{subsection}{2}{\z@}%
  {1.0ex plus 0.4ex minus 0.2ex}%
  {0.5ex plus 0.2ex minus 0ex}%
  {\normalfont\normalsize\itshape}}%
\makeatother

\makeatletter
\def\@IEEENORMtitlevspace{2.15\baselineskip}
\def\@IEEEMINtitlevspace{1.65\baselineskip}
\makeatother

\begin{document}

\IEEEaftertitletext{\vspace{-0.35\baselineskip}}

\title{\fontsize{18}{20}\selectfont Quantum Restricted Boltzmann Machine for Fast Unit Commitment}

\author{Yuji Cao,~\IEEEmembership{Graduate Student Member,~IEEE,}
Yue Chen,~\IEEEmembership{Senior Member,~IEEE}
\thanks{This work was supported by the Hong Kong Innovation and Technology Commission Innovation and Technology Fund under Grant ITS/335/24.}
\thanks{Y. Cao and Y. Chen are with the Department of Mechanical and Automation Engineering, The Chinese University of Hong Kong, Hong Kong, China. (e-mail: \{yjcao, yuechen\}@mae.cuhk.edu.hk). Corresponding author: Y. Chen.}}


\maketitle

\begin{abstract}
Network-constrained unit commitment (UC) problems must be solved repeatedly to ensure power system reliability, yet face computational challenges from the growing system scale and fast response requirements. This letter proposes a quantum restricted Boltzmann machine (qRBM) method to accelerate UC solution by predicting solution patterns via quantum-hardware-compatible energy scoring and recovering optimal solutions through affine mapping. This method reduces the computation time for UC  by 99.96\%. Compared to existing quantum methods, it cuts the number of qubits needed by 276 times.
\end{abstract}


\section{Introduction}

In practice, unit commitment (UC) must be solved repeatedly to accommodate changing renewable forecasts, but the underlying mixed-integer optimization problem remains computationally challenging due to combinatorial commitment decisions and growing network scale. Quantum computing shows great potential for combinatorial optimization and offers a promising solution to this challenge~\cite{Nikmehr2022, 11449549}.
However, most existing approaches reformulate UC as quadratic unconstrained binary optimization (QUBO) via binary expansion, which demands qubit resources far exceeding current quantum hardware capability, even for small instances~\cite{Nikmehr2022}.

To enable large-scale UC solution on real quantum hardware, this letter proposes a quantum restricted Boltzmann machine (qRBM) method that first predicts UC solution patterns (i.e., commitment decisions and active constraint sets) via quantum-hardware-compatible energy scoring, and then recovers optimal solutions through affine mapping. By encoding learned pattern partitions rather than the full problem formulation, this method significantly reduces qubit requirements and improves scalability on near-term quantum devices. Furthermore, the trained qRBM can be deployed directly on real Ising-type quantum hardware and is practical.

\section{Unit Commitment Problem Reformulation}
\label{sec:formulation}

We consider a network-constrained unit commitment (UC) problem over $T$ periods, parameterized by renewable forecasts. The problem is a mixed-integer quadratic program (MIQP):
\begin{subequations}
\label{eq:miqp}
\begin{align}
\min_{z,u,v, P}~&
\sum \nolimits_{t=1}^{T}\sum \nolimits_{i=1}^{N}
\bigl(F_i^C(P_{it})+F_i^U(z_{it},u_{it},v_{it})\bigr), \label{eq:miqp-obj}\\
\text{s.t.}~&
(z,u,v)\in\mathcal{Y},~ z, u, v \in \{0,1\}^{NT}, ~ P\in\Omega(z),
\label{eq:miqp-ucset}\\
&
\sum \nolimits_i P_{it}+\sum \nolimits_k\theta_{kt}=\sum \nolimits_j D_{jt},
\  \forall t, \label{eq:miqp-balance}\\
& -f^{\max}\le \Pi_g P_t+\Pi_\theta\theta_t-\Pi_dD_t\le f^{\max},
\  \forall t. \label{eq:miqp-flow}
\end{align}
\end{subequations}
Here, $P_{it}$ is the power output of generator $i$ at time $t$; $z_{it}$, $u_{it}$, and $v_{it}$ are commitment, startup, and shutdown binary decisions; $\theta_{kt}$ and $D_{jt}$ are renewable output and load demand; and $\Pi_g$, $\Pi_\theta$, and $\Pi_d$ are power transfer distribution factors. The quadratic functions $F_i^C$ and $F_i^U$ denote generation and commitment-related costs. $\mathcal{Y}$ is the feasible set of UC decisions, consisting of startup/shutdown logic and minimum up/down times, while $\Omega(z)$ collects the generation capacity and ramping limits. Constraints \eqref{eq:miqp-balance} and \eqref{eq:miqp-flow} are the power balance condition and power network constraints, respectively.

Each update of the renewable forecast  $\theta$ necessitates re-solving the UC problem \eqref{eq:miqp}.
However, due to the large number of binary variables, solving~\eqref{eq:miqp} is time-consuming. In the following, we propose a reformulation of \eqref{eq:miqp}, based on which the quantum restricted Boltzmann machine method is developed to accelerate the solution process.
Let $y=(z_{it},u_{it},v_{it},\forall i,\forall t)$ collect the binary variables and $x=(P_{it},\forall i, \forall t)$ collect the continuous variables. For a given $\theta$, if a binary vector $y$ is known, the commitment costs become constant and~\eqref{eq:miqp} reduces to a parametric QP in~$x$:
\begin{equation}
\min_x \; \tfrac{1}{2} x^\top Q x + c_x^\top x, \quad
\text{s.t.} \; Ax \le (e - Ey) + B\theta,
\label{eq:qp}
\end{equation}
where $Q \succ 0$ is block-diagonal with entries $2a_i$, $c_x$ collects the linear cost coefficients~$b_i$, and $A$, $e$, $B$, $E$ collect the constraint coefficients.

By multi-parametric programming~\cite{pistikopoulos2020multi}, for a fixed binary vector $y$, the optimal solution $x^*(\theta)$ of~\eqref{eq:qp} is piecewise affine in~$\theta$. Based on this observation, the feasible region of the uncertainty parameter $\theta$ can be partitioned into a set of \textit{critical regions} (CRs), such that for all $\theta$ within the same CR, both the optimal binary variables and the active constraint set of the resulting problem \eqref{eq:qp} remain unchanged; the $\ell$-th critical region is defined by
\begin{equation}
\mathrm{CR}_\ell
:=
\bigl\{\theta:\,
y^*(\theta)=y_\ell,\;
\mathcal{A}^*(\theta)=\mathcal{A}_\ell
\bigr\},
\label{eq:cr}
\end{equation}
where $y_\ell$ is the optimal binary vector and $\mathcal{A}_\ell$ is the active constraint set of problem~\eqref{eq:qp}. We also call $(y_\ell,\mathcal{A}_\ell)$ a UC pattern. Within each critical region, the continuous dispatch $x^*$ is affine in~$\theta$ and can be recovered without re-solving~\eqref{eq:qp}. Specifically, substituting $y=y_\ell$ into \eqref{eq:qp} and deriving its Karush–Kuhn–Tucker (KKT) conditions yields:
\begin{equation}
\underbrace{\begin{bmatrix} Q & A_{\mathcal{A}_{\ell}}^\top \\ A_{\mathcal{A}_{\ell}} & 0
\end{bmatrix}}_{K_\ell}
\begin{bmatrix} x^* \\ \lambda^* \end{bmatrix}
= \underbrace{\begin{bmatrix} -c_x \\ e_{\mathcal{A}_{\ell}} 
    - E_{\mathcal{A}_{\ell}} y_{\ell} \end{bmatrix}}_{d_\ell} + \underbrace{\begin{bmatrix}
        0 \\ B_{\mathcal{A}_{\ell}}
    \end{bmatrix}}_{R_\ell}\theta,
\label{eq:kkt}
\end{equation}
where $A_{\mathcal{A}_{\ell}}$, $e_{\mathcal{A}_{\ell}}$, $B_{\mathcal{A}_{\ell}}$, $E_{\mathcal{A}_{\ell}}$ are the active-set rows of $A$, $e$, $B$, $E$ and $\lambda^*$ contains the optimal dual variables. Denote the left-hand side coefficient matrix in~\eqref{eq:kkt} by~$K_\ell$ and the right-hand side coefficients by $d_\ell$ and $R_\ell$, respectively. Inverting~$K_\ell$ and extracting the first~$n_x$ rows gives the dispatch solution:
\begin{equation}
x^*(\theta) = M_\ell\, \theta + m_\ell, \quad y^*(\theta) = y_\ell, \quad
\theta \in \mathrm{CR}_\ell,
\label{eq:affine}
\end{equation}
with $M_\ell = [K_\ell^{-1} R_\ell]_x$ and $m_\ell = [K_\ell^{-1} d_\ell]_x$, where $[\cdot]_x$ selects the first~$n_x$ rows and $\mathrm{CR}_\ell$ is the $\ell$-th critical region. Therefore, once the critical region containing $\theta$ is identified, $x^*(\theta)$ and $y^*(\theta)$ can be recovered directly from~\eqref{eq:affine}. The UC problem~\eqref{eq:miqp} is thereby recast as a critical region prediction problem.

To construct the candidate CRs, the UC problem~\eqref{eq:miqp} is solved for each sampled renewable forecast. Each optimal solution yields a pair $(y^*,\mathcal{A}^*)$, and all distinct pairs are enumerated as $L$ critical regions. For each $\mathrm{CR}_\ell$, the matrices $M_\ell$ and $m_\ell$ in~\eqref{eq:affine} are precomputed from~\eqref{eq:kkt}.

\section{Quantum Boltzmann Machine for UC}
\label{sec:method}

In this section, we develop a qRBM for predicting the critical region $\mathrm{CR}_{\ell}$ for a given $\theta$. Once the optimal CR appears among the top-ranked candidates,~\eqref{eq:affine} recovers the exact optimal UC solution.

Specifically, the qRBM serves as a pattern scorer: for each candidate CR, the model assigns a free-energy score conditioned on the input $\theta$. The lower the free energy, the more likely $\theta$ belongs to the candidate CR. The qRBM consists of two layers: 1) a visible layer encoding the input and output information, and 2) a hidden layer capturing latent dependence between them. In this letter, we use the qRBM architecture in~\cite{Demidik2025}, where the visible layer stores classical binary information and the hidden layer is represented by quantum units. We follow three steps:

\textbf{\textit{Step 1: State Preparation}}. Before applying the qRBM, we encode each input parameter (renewable generation) and candidate CR as one visible state. The renewable input is encoded as $\xi=\operatorname{enc}(\theta)\in\{0,1\}^{n_\xi}$ by binary expansion. The pre-collected $L$ candidate CRs are represented by one-hot codes $c_\ell\in\{0,1\}^{L}$; for $\mathrm{CR}_{\ell}$, $c_{\ell,\ell}=1$ and $c_{\ell,k}=0,\forall k \ne\ell$.  Combining them, for $\mathrm{CR}_\ell$, the corresponding visible state is
\begin{equation}
  r_\ell(\xi)=[\xi,c_\ell]\in\{0,1\}^{n_{\text{vis}}},
  \quad n_{\text{vis}}=n_\xi+L .
  \label{eq:visible}
  \end{equation}
After encoding, each visible bit is mapped to an Ising spin, i.e., $s_{\ell,n}=1-2r_{\ell, n}\in\{-1,+1\},\forall n=1,\cdots,n_{\text{vis}}$.

\textbf{\textit{Step 2: qRBM Training}}. Given a visible state $r$, let $\psi$ denote all trainable qRBM parameters. The qRBM free energy $\mathcal{G}_\psi(r)$ is defined as
\begin{equation}
\mathcal{G}_\psi(r)
=
\sum \nolimits_{i=1}^{n_{\text{vis}}}s_i (r_i) \alpha_i
-
\sum \nolimits_{j=1}^{n_{\text{hid}}}
\log\cosh\bigl(\|\Phi_j(r)\|_2\bigr),
\label{eq:qrbm_free_energy}
\end{equation}
where $s_i(r_i)=1-2r_i$ as above, $\alpha_i$ is a trainable visible bias, $n_{\mathrm{hid}}$ is the number of quantum hidden units, and $\Phi_j(r)$ measures how strongly the visible state $r$ activates the $j$-th hidden unit. The first term in~\eqref{eq:qrbm_free_energy} scores the visible bits of $r$; the second term aggregates the hidden-unit contribution induced by $r$. For a given renewable input, a more probable CR receives a lower free energy $\mathcal{G}_\psi(r)$ and is ranked higher.

To train the qRBM, we first convert the free energy score into the unnormalized probability below: 
\begin{equation}
  \tilde p_\psi(r)=\exp[-\mathcal{G}_\psi(r)].
  \label{eq:qrbm_prob}
\end{equation}
Thus, lower free energy gives larger probability. For a given input~$\xi$, the conditional probability of $\mathrm{CR}_\ell$ is defined over valid one-hot labels as
\begin{align}
p_\psi(\ell\mid \xi)
&=
\frac{\tilde p_\psi(r_\ell(\xi))}
{\sum_{i=1}^{L}\tilde p_\psi(r_i(\xi))},\\
\mathcal{L}(\psi)
&=
-\sum \nolimits_{(\xi,\ell^*)\in\mathcal D}
\log p_\psi(\ell^* \mid \xi).
\label{eq:loss}
\end{align}
The qRBM is trained on classical hardware by minimizing $\mathcal{L}(\psi)$, where $\mathcal D$ is the training set and $\ell^*$ is the CR label corresponding to the MIQP optimal solution. Training adjusts the parameters $\psi$ so that the correct CR receives high probability for each input parameter (renewable power) $\theta$ (or $\xi$).

\textbf{\textit{Step 3: UC Solution}}. The trained qRBM is used to accelerate the UC solution process. Specifically, given a renewable input~$\theta$, we encode it as $\xi=\operatorname{enc}(\theta)$, form the visible states $r_\ell(\xi),\forall \ell=1,\cdots,L$, and map them to Ising spin vectors $s_{\ell}(\xi)=1-2r_{\ell}(\xi) \in \{-1,+1\}^{n_{\text{vis}}}$. The trained parameters~$\psi$ define the corresponding Ising Hamiltonian on the quantum computer. The quantum computer samples low-energy states, which correspond to candidate CRs with large $p_\psi(\ell\mid\xi)$. We select the $n_{\mathrm{cand}}$ most likely CRs and then choose the lowest-cost feasible one:
%
\begin{equation}
  \begin{aligned}
  \widehat{\mathcal S}(\xi)
  &= \operatorname*{argmax}_{n_{\mathrm{cand}},\,\ell\in\{1,\ldots,L\}}
     p_\psi(\ell\mid \xi), \\
  \hat\ell
  &= \arg\min_{\ell \in \widehat{\mathcal S}(\xi)}
     \bigl\{J(\ell,\theta): M_\ell\theta+m_\ell
     \text{ is feasible}\bigr\}.
  \end{aligned}
  \label{eq:candidate_selection}
\end{equation}
where $\operatorname*{argmax}_{n_{\mathrm{cand}}}$ returns the indices of the $n_{\mathrm{cand}}$ largest values, and $J(\ell,\theta)$ is the objective value after applying the affine solution mapping~\eqref{eq:affine}, given $(\ell,\theta)$.
Each candidate is recovered by~\eqref{eq:affine} and checked in~\eqref{eq:miqp}. If a feasible candidate exists in $\widehat{\mathcal S}(\xi)$,~\eqref{eq:candidate_selection} returns the lowest-cost one. If the UC pattern of~\eqref{eq:miqp} lies in $\widehat{\mathcal S}(\xi)$, the returned solution is optimal.

\begin{table*}[thbp]
  \caption{Network-constrained UC Performance across Different Systems (Unit: \%)}
  \label{tab:perf}
  \centering
  \scriptsize
  \setlength{\tabcolsep}{6pt}
  \begin{tabular}{@{} l cc cc cc @{}}
  \toprule
  \multirow{2}{*}[-0.5ex]{\textbf{Method}}
    & \multicolumn{2}{c}{\textbf{9-bus}}
    & \multicolumn{2}{c}{\textbf{57-bus}}
    & \multicolumn{2}{c@{}}{\textbf{118-bus}} \\
  \cmidrule(lr){2-3}\cmidrule(lr){4-5}\cmidrule(l){6-7}
  & Feasibility Rate & Opt. Rate
    & Feasibility Rate & Opt. Rate
    & Feasibility Rate & Opt. Rate \\
  \midrule
M1: ANN direct prediction
  & 27.37 $\pm$ 2.33 & 20.36 $\pm$ 0.83 & 0.00 $\pm$ 0.00 & 0.00 $\pm$ 0.00 & 28.14 $\pm$ 1.75 & 23.94 $\pm$ 1.48 \\
M2: ANN direct prediction with simple repair
  & 74.11 $\pm$ 1.87 & 36.95 $\pm$ 2.92 & 88.14 $\pm$ 1.40 & 14.94 $\pm$ 0.26 & 44.26 $\pm$ 1.16 & 30.50 $\pm$ 1.74 \\
M3: qRBM direct prediction with simple repair
  & 71.69 $\pm$ 1.33 & 42.58 $\pm$ 1.45 & 97.88 $\pm$ 0.41 & 15.34 $\pm$ 0.15 & 51.68 $\pm$ 1.74 & 45.72 $\pm$ 0.96 \\
M4: VQC with reformulation
  & 93.54 $\pm$ 1.33 & 45.08 $\pm$ 1.65 & 99.86 $\pm$ 0.15 & 58.94 $\pm$ 1.25 & 88.56 $\pm$ 0.82 & 83.20 $\pm$ 1.26 \\
\textbf{Ours: qRBM with reformulation}
  & \textbf{100.00} $\pm$ \textbf{0.00} & \textbf{97.66} $\pm$ \textbf{0.50} & \textbf{100.00} $\pm$ \textbf{0.00} & \textbf{96.08} $\pm$ \textbf{0.15} & \textbf{97.38} $\pm$ \textbf{0.33} & \textbf{93.32} $\pm$ \textbf{0.35} \\
  \bottomrule
  \end{tabular}
\end{table*}

\section{Case Studies}

The proposed method is validated on modified IEEE 9-bus, 57-bus, and 118-bus systems with a 24-hour horizon. For each system, we solve UC problems under sampled renewable scenarios and construct the datasets. The sample sizes for 9-bus, 57-bus, and 118-bus are 3,840, 80,000, and 97,412, respectively, with an 80\%/20\% training-test split.

\subsection{Overall UC Performance}

We compare the proposed method with four benchmarks: 1)~\textbf{M1}: artificial neural network (ANN) with direct prediction \cite{shokry2017mixed}. For binary variables, it predicts directly. For continuous variables, k-means clustering and linear regressions are applied to approximate the division of critical regions and the mapping between continuous variables
and optimal solutions, respectively. 2)~\textbf{M2}: ANN with direct prediction and simple repair\footnotemark; 3)~\textbf{M3}: qRBM with direct prediction and simple repair; 4)~\textbf{M4}: variational quantum circuit (VQC) with the proposed reformulation. 

\footnotetext{Simple repair: Unit commitment decisions are projected to the feasible region to meet constraints \eqref{eq:miqp-ucset}-\eqref{eq:miqp-flow}.}

\begin{figure}[thbp]
  \centering
  \includegraphics[width=1.00\columnwidth]{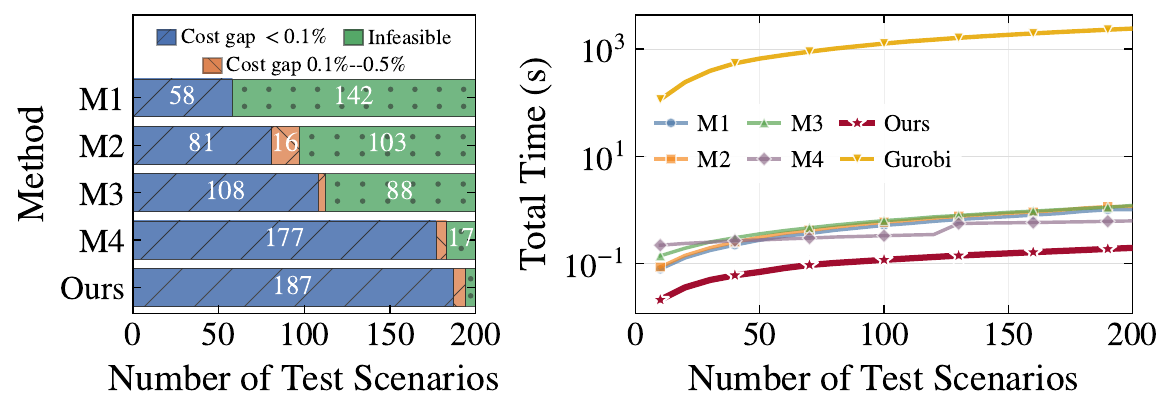}
  \caption{Cost gap and total time comparison on the 118-bus system.}
  \label{fig:comparison}
\end{figure}

Table~\ref{tab:perf} reports the mean feasibility and optimality rates on the three systems. Feasibility and optimality rates are the share of test scenarios with a feasible recovered solution and an objective matching the optimum of~\eqref{eq:miqp}, respectively. Direct methods (M1--M3) have poor performance, especially in optimality, whereas reformulation improves both M4 and the proposed method. Compared to VQC, the proposed method reaches higher feasibility and greatly improves the optimality on all three systems. Fig.~\ref{fig:comparison} further zooms in on 200 scenarios of the 118-bus system. The proposed method keeps the cost gap below $0.1\%$ in most cases (187 / 200), and its cumulative time stays well below the baselines and orders of magnitude below Gurobi.


\subsection{qRBM Learned Patterns}
  


\begin{figure}[thbp]
\centering
\includegraphics[width=1.00\columnwidth]{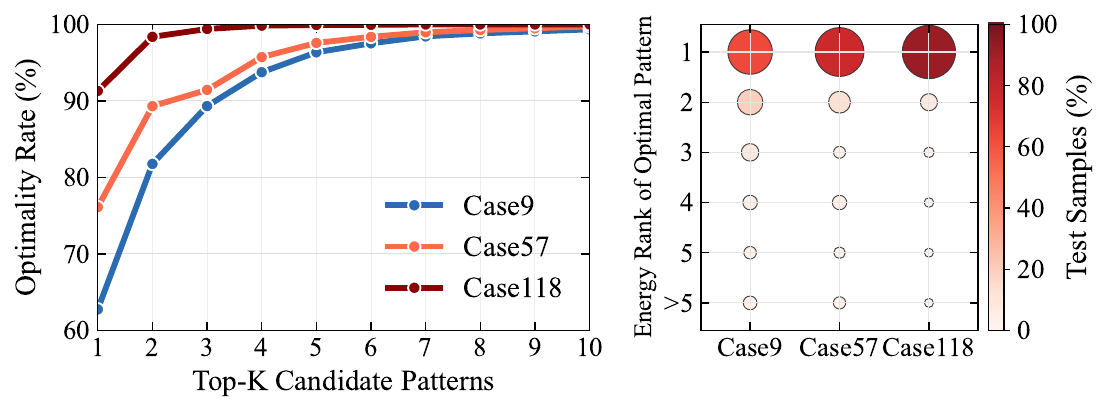}
\caption{Top-$n_{\text{cand}}$ optimality and optimal-CR energy rank by learned qRBM.}
\label{fig:topk}
\end{figure}

Next, we assess whether the qRBM assigns low energy to the CR corresponding to the optimal solution. For each test sample, candidate CRs are ranked by qRBM energy, and we check whether the MIQP-optimal CR appears among the top-\(n_{\text{cand}}\) candidates.
Fig.~\ref{fig:topk} (left) shows that the optimality rate increases rapidly with $n_{\text{cand}}$. Notably, the top-5 low-energy candidates already include an MIQP-optimal CR for over 95\% of test samples. Fig.~\ref{fig:topk} (right) reports the energy rank at which the MIQP-optimal CR appears. Bubble size and color indicate the percentage of test samples at each rank. Large bubbles on rank 1 show that the optimal CRs are concentrated among the top-ranked low-energy candidates.

\subsection{Performance on Real Quantum Hardware}

\begin{table}[hbtp]
  \caption{Computation Time and Qubit Requirement Comparison}
  \label{tab:scale}
  \centering
  \scriptsize
  \setlength{\tabcolsep}{1.6pt}
  \begin{tabular*}{\columnwidth}{@{\extracolsep{\fill}} lcccccc @{}}
  \toprule
  System
    & \makecell{Qubits\\(QUBO)}
    & \textbf{\makecell{Qubits\\(Ours)}}
    & \textcolor{black}{\makecell{Qubit\\Efficiency}}
    & \makecell{Gurobi\\Time (ms)}
    & \textbf{\makecell{Ours\\Time (ms)}}
    & \textcolor{black}{\makecell{Time\\Speedup}} \\
  \toprule
  \textcolor{black}{9-bus}    & 4,513    & \textbf{46}  & \textcolor{black}{98.1$\times$}  & 670   & \textbf{0.778} & \textcolor{black}{861$\times$} \\
  \textcolor{black}{57-bus}   & 33,185   & \textbf{113} & \textcolor{black}{293.7$\times$} & 1,960  & \textbf{0.197} & \textcolor{black}{9,949$\times$} \\
  \textcolor{black}{118-bus}  & 90,433   & \textbf{381} & \textcolor{black}{237.4$\times$} & 11,300 & \textbf{1.121} & \textcolor{black}{10,080$\times$} \\
  \textcolor{black}{1354-bus} & 212,457  & \textbf{448} & \textcolor{black}{474.2$\times$} & 13,820 & \textbf{2.739} & \textcolor{black}{5,046$\times$} \\
  \bottomrule
  \end{tabular*}
\end{table}

\begin{figure}[thbp]
\centering
\includegraphics[width=1.00\columnwidth]{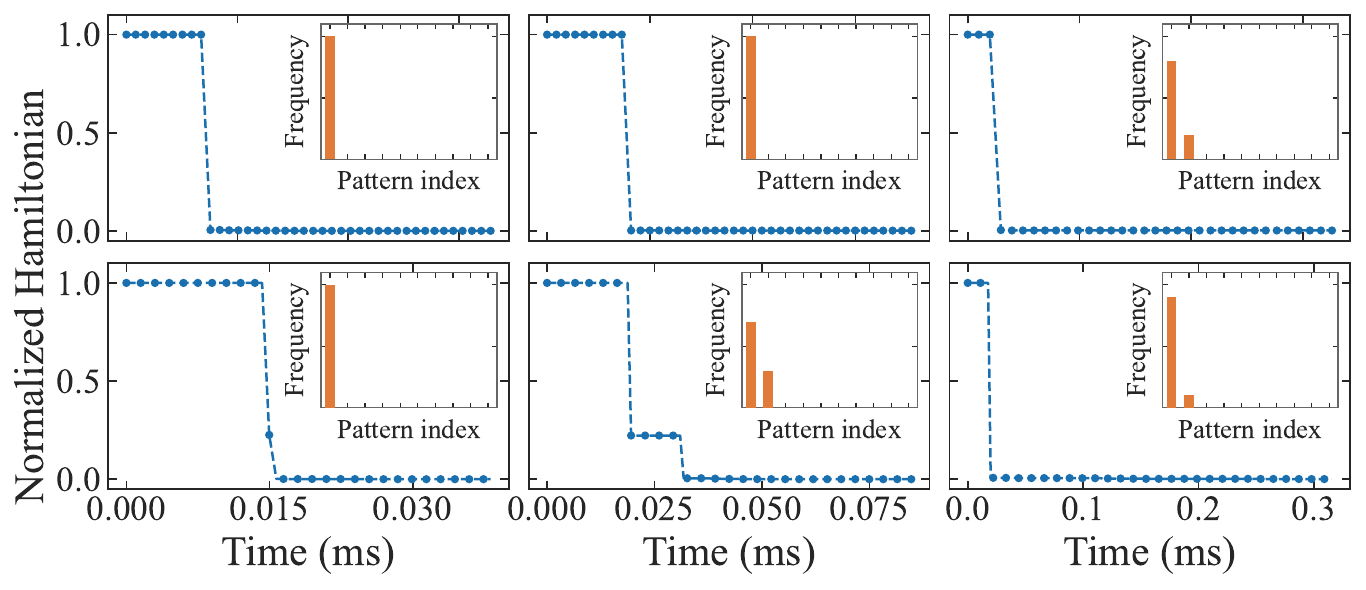}
\caption{Ising-Hamiltonian evolution on the QBoson coherent Ising machine.}
\label{fig:qpu}
\end{figure}

We further test the proposed method on a modified IEEE 1354-bus system using a real quantum computer CPQC-X, which is a coherent Ising machine (CIM) from QBoson. 
Fig.~\ref{fig:qpu} shows the Ising-Hamiltonian energy evolutions of six representative test scenarios. In all tested scenarios, the CIM rapidly converges to a low-energy state and returns 10 candidates, from which the optimal CR can be found. This result validates the practicability of applying the proposed method on real quantum hardware.

Table~\ref{tab:scale} compares qubit requirement and computation time against a direct QUBO encoding and Gurobi solver. From the 9-bus to the 1354-bus systems, a QUBO formulation needs $4{,}513$--$212{,}457$ qubits, whereas the proposed method uses only $46$--$448$ qubits, because the qRBM hardware size follows the learned CR representation rather than the full UC formulation. The resulting encoding efficiency is $98.1\times$--$474.2\times$ (mean $276\times$). The computation time of our method (CIM plus affine recovery) remains within 0.778--2.739 ms across all systems, compared with 670--13{,}820 ms for Gurobi, corresponding to a $861\times$--$10{,}080\times$ speedup and reducing computation time by 99.96\% on average.

\section{Conclusion}

This letter proposes a novel qRBM method with solution pattern prediction for UC problems, enabling large-scale solution with far fewer qubits and better optimality than current quantum formulations.

\bibliographystyle{IEEEtran}
\bibliography{mybib}

\end{document}